\documentclass{article}
\usepackage{amsmath,amssymb,amsthm}
\usepackage{graphics}

\newtheorem{theorem}{Theorem}

\title {CONSTRUCTION OF THE VORONOI DIAGRAM AND SECONDARY POLYTOPE}

\author {Oleg R. Musin \thanks{The research was partially supported by the
RFFI grant 97-01-00174 and by the CTIAC ESPRIT Project No. 21042}
\\{ } \\
Moscow State University, Moscow, 119899, Russia\\
e-mail: omusin@mail.ru}

\begin{document}
\date{}
\maketitle


\centerline{\bf ABSTRACT}

\medskip
A set $S$ of $n$ points in general position in ${\bf R}^d$ defines the unique
Voronoi diagram of $S$. Its dual tessellation is the Delaunay triangulation
(DT) of $S.$ In this paper we consider the parabolic functional on the set of triangulations
of $S$ and prove that it attains its minimum at DT in all dimensions. The Delaunay 
triangulation of $S$ is corresponding to a vertex of the secondary
polytope of $S$. We proposed an algorithm for DT's construction, where the
parabolic functional and the secondary polytope are used.
Finally, we considered a discrete analog of the Dirichlet functional.
DT is optimal for this functional only in two dimensions.
\medskip


\section {Introduction}
Some of the most well-known names in Computational Geometry are those
of two prominent Russian mathematicians:  Georgy F. Voronoi
and Boris N. Delaunay.  Their considerable contribution to the
Number Theory and Geometry  is well known to the specialists in these fields.
Surprisingly, their names (their works remained unreaded and later
re-discovered) became the most popular not among "pure" mathematicians,
but among the researchers who used geometric applications. Such terms as
"Voronoi diagram" and "Delaunay triangulation" are very important not only
for Computational Geometry, but also for Geometric Modeling, Image Processing,
CAD, Geographic Information System  etc.

The Voronoi diagram is generated by a set of $n$ points $S=\{x_1,...,x_n\}$
in ${\bf R}^d$. The {\em Voronoi diagram} is the partition of the ${\bf R}^d$
into $n$ convex cells, the Voronoi cells $V_i$, where each $V_i$ contains all
points of the ${\bf R}^d$ closer to $x_i$ than to any other point:
$$ V_i = \{x|\forall j \ne i, d(x,x_i) \le d(x,x_j)\},$$ where $d(x,y)$ is
the Euclidean distance between $x$ and $y$.

This concept has independently appeared in various fields of science.
The earliest significant use of Voronoi diagrams seems to have occurred in
the work of Gauss, Dirichlet and Voronoi in their investigations on the
reducibility of positive definite quadratic forms.
Different names particular to the respective field have been used, such as
{\em medial axis transform} in biology and physiology, {\em Wigner-Seitz
zones} in chemistry and physics, {\em domain of action} in crystallography,
and {\em Thiessen polygons} in meteorology and geography [1].

For generic set of $n$ points $S$ in ${\bf R}^d$ the straight-line dual of the
Voronoi diagram is triangulation of $S$, called the {\em Delaunay
triangulation} and denoted by $DT(S).$  The $DT(S)$ is triangulation of the
convex hull of S in ${\bf R}^d$ and set of vertices of $DT(S)$ is $S$.

Voronoi [21] was the first to consider the dual structure of the Voronoi
diagram, where any two points of $S$ are connected whose regions have a
boundary in common. Later Delaunay [6] obtained the same by defining that
two points of $S$ are connected if and only if they lie on a sphere whose
interior contains no points of $S$.  After him, the dual of the Voronoi
diagram has been denoted {\em Delaunay tessellation} or {\em Delaunay
triangulation}.

Voronoi diagrams and Delaunay triangulations are used in numerous
applications. It is widely used in plane and 3D case.
A natural question may arise:  why those structures are better than the
others. Usually the advantages of planar Delaunay triangulation are
rationalized by the max-min angle criterion and other properties
[1,2,3,6,16,17,19].

The max-min angle criterion requires that the diagonal of
every convex quadrilateral occurring in the triangulation "should be well
chosen" [19], in the sense that replacement of the chosen diagonal by the
alternative one must not increase the minimum of the six angles in the two
triangles making up the quadrilateral. Thus the Delaunay triangulation of a
planar point set maximizes the minimum angle in any triangle. More
specifically, the sequence of triangle angles, sorted from sharpest to least
sharp, is lexicographically maximized over all such sequences constructed from
triangulation of $S$.

In the papers [13--15] we defined several functionals on the set of all
triangulations of $S$ in ${\bf R}^2$ attaining global minimum on the
Delaunay triangulation.

The "mean radius" functional is the mean of circumradii of
triangles for planar triangulations. Let $t$ be a triangulation of
$S$ in the plane. Assume  that each triangle $\Delta_i$ of this triangulation
is related to the radius $R_i$ of its circumcircle. Thus every triangulation
$t$ is related to the set $\{R_{\Delta_1},...,R_{\Delta_k}\}$ of circumradii
of triangles $\Delta_i \in t$. The numbers of triangles for
any two triangulations of $S$ are equal, so it is possible to compare sets
of radii for different triangulations. In particular, it is possible to
compare sums of radii: $\sum R_{\Delta_i}$ or power sums:
$\sum R_{\Delta_i}^a, a>0$. It seems that  triangulation having
minimal sum of radii is "better", because all its triangles in "average"
are nearer to the regular triangles. The functional
$R(t,a) = \sum R_{\Delta_i}^a, a>0$ attains its minimum iff $t$ is
Delaunay triangulation [15].

The harmonic index of triangulation has its origin in the theory of the so
called "harmonic maps". For polygon $P$ its {\em harmonic index}
$$hrm(P) = \sum a_i^2/S(P),$$ where $a_1,\ldots, a_n$ are the
lengths of sides of $P$ and $S(P)$ is its area. This index is the same
for similar polygons. It is easy to prove that harmonic index for triangles
achieves its minimum iff triangle is equiangular.
For planar triangulation $t$ of a set $S$ let denote by $hrm(t)$
(harmonic index of triangulation $t$) the sum of $hrm$ of its triangles:
$$hrm(t)=\sum\limits_{\Delta_i \in t}hrm(\Delta_i)$$. Harmonic index
$hrm(t)$ of triangulation $t$ of $S$ achieves its minimum iff $t$ is the
Delaunay triangulation of $S$.
The harmonic functional for triangle attains its minimum if this triangle
is equiangular. Usually, a triangulation is regarded as "good" for different
purposes if its triangles are nearly equiangular. The harmonic index of
triangulation $t$ achieves its global minimum if $t$ is a part of a regular
triangular lattice. In some sense, this result shows that the Delaunay
triangulation is as close as possible to  equiangular triangulation.

One of the most popular algorithm for constructing planar Voronoi diagram
-- Delaunay triangulation is so called flipping algorithm [1,8,10].
(If $t$ is planar triangulation of $S$ and $AC$ is an internal edge of $t$,
the two triangles $ACB$ and $ACD$ of $t$ incident with $AC$ in
triangulation. If quadrilateral $ABCD$ is convex a new triangulation of $S$
may be obtain by removing the edge $AC$ and inserting the edge $BD$. This
operation is called the {\em flip} or {\em edge--flip}).
The proof that sequence of flips does not cycle in the algorithm easily
follows from consideration of the mean-radius or harmonic functional.
Indeed, using flipping algorithm after each flip the functional decreases
until Delaunay triangulation is reached.

Among the various proposed methods for constructing Voronoi diagrams in
${\bf R}^d$, incremental insertion of points is most intuitive and easy
implement. Joe [9], Rajan [16], and Edelsbrunner and Shah [8]
generalized of planar incremental and flipping algorithms for higher
dimensions.

This paper includes algorithm for construction of the Voronoi diagram.
First, in section 1 we consider a so called parabolic functional on a set
of all triangulations of $S$ and prove that this functional attains its
minimum on the Delaunay triangulation of $S$ in all dimensions.
The secondary polytope is the original one due Gel'fand, Kapranov,
and Zelevinsky [4,5]. They introduced the secondary polytope
$Q$ of $S \subset {\bf R}^d$ that is a convex polytope in ${\bf R}^n, n=|S|$,
and where the vertices of $Q$ are in one-to-one correspondence with the
regular triangulations of $S.$ In section 2 we show that parabolic functional
is a linear function on ${\bf R}^n$, and Delaunay triangulation
of $S$ is a vertex of $Q$ that gives minimum for parabolic function.
The incremental algorithm for construction of the Voronoi diagram i.e.
Delaunay triangulation is given in section 3. The main idea of this algorithm
is to construct a sequence of a regular triangulations (vertices of $Q$) that
decrease parabolic functional, and the last triangulation in this sequence is
DT. Finally, in section 4 we consider a discrete analoque of the Dirichlet
functional on a set of all triangulations of $S.$ For $d=2$ this functional
achieves its minimum for DT. We consider case $n=d+2$ and prove that optimal
triangulation depends on points $S$ configuration only.
If $d>2$ then Delaunay triangulation is not
optimal for Dirichlet functional. Thus the problem to finding "good"
triangulations for this functional in higher dimensions is opened and more
detailed consideration is necessary.

\section {Optimality of Delaunay triangulations for the parabolic functional}

    Throughout this paper $S=\{x_1,...,x_n\}$ denotes a set of $n$ points in
general position in ${\bf R}^d$. A {\em triangulation} of $S$ is a
triangulation of the polytope $P=CH(S)$ (convex hull of $S$) with
vertices in $S.$

Let $t$ be a triangulation of the set $S$ in ${\bf R}^d$,
$\Delta_i$ denotes the $i$-th $d$-simplex of $t$ and
$x_{ij} \in {\bf R}^d, \, j = 0, 1, \ldots, d$ are its vertices .
Let $$Vr(t) = \sum\limits_i (x_{i0}^2+...+ x_{id}^2) vol(\Delta_i),$$
where $vol(\Delta_i)$ is volume (area for $d=2$) of the simplex $\Delta_i$.
We call functional $Vr$ {\em parabolic} (or {\em Voronoi}).

The parabolic functional induces an order on triangulations of the set
$S$ by the rule: $t_1 > t_2$ iff $Vr(t_1) > Vr(t_2)$. The value of $Vr$
depends on the choice of the origin. If we move the origin to $x_0$
then this order does not change i.e.

\medskip
$Vr(t_1) > Vr(t_2)$ iff $Vr(t_1,x_0) > Vr(t_2,x_0)$.
\medskip

The main result for Voronoi functional $Vr$ is the following:

\begin{theorem} The parabolic (Voronoi) functional $Vr(t)$ achieves its
minimum if and only if $t$ is the Delaunay triangulation
\end{theorem}

A simple proof of these theorems follows from paroboloid construction of
the DT found by Voronoi [21] and rediscovered only in 1979 for a sphere
(K.Brown), and later also for a paraboloid.

Let us consider an arbitrary triangulation $t$ of the set $S$ and "lift" it
onto the paraboloid in ${\bf R}^{d+1}$, i.e.
let us build a polyhedral surface in ${\bf R}^{d+1}$ connecting corresponding
vertices on the paraboloid. Note that the functional $Vr$  up to a constant
equals to the volume of the solid body below this surface. Thus, the minimum
of $Vr$ is attained on the Delaunay triangulation.

Let us consider another functional for triangulations:
$$ C2(t) = \sum\limits_i ||c(\Delta_i)||^2 vol(\Delta_i), $$
where $c(\Delta_i)$ is the center (barycenter) of the $\Delta_i$,
$c(\Delta_i) = \sum_j x_{ij}/(d+1)$ and $||\cdot||$ is Euclidean norm.

By direct calculation (it is sufficient to check the formula on a simplex)
it is possible to prove that
$$(d+1)^2 C2(t) + Vr(t) = (d+1)(d+2)\int\limits_{CH(S)}||x||^2dx,$$
where $CH(S)$ is convex hull of set $S$ in ${\bf R}^d$.

From Theorem 1 and this formula directly follows that:

\begin{theorem} The functional C2(t) on triangulations of the set $S$
achieves its maximum if and only if $t$ is the Delaunay triangulation.
\end{theorem}

\section {The secondary polytope, regular triangulations, and flips}

The following analytic description of the secondary polytope is
the original one due Gel'fand, Kapranov, and Zelevinsky [4,5]. We include it
here for completeness.
They introduced the secondary polytope $Q = \sum(S)$ of an affine point
configuration $S$, where the vertices of $\sum(S)$ are in one-to-one
correspondence with the regular triangulations of the "primary polytope"
$P = CH(S)$ - convex hull of $S$.

Let $S = \{x_i\}$ be a set of $n$ points in ${\bf R}^d$ and $t$ is a
triangulation of $S$. We are correspond to triangulation $t$ vector
$p(t) = (p_1, p_2, \ldots, p_n)$ in ${\bf R}^n,$ where
$p_i = \sum vol(\Delta_{ij}, \Delta_{ij}$ denotes the $j$-th $d$-simplex
of $t$ that incident to vertex $x_i,$ and $vol(\Delta)$ denotes volume of
the $d$-simplex i.e. $p_i$ is volume of star of $i$-th vertex of
triangulation $t$.

We have $d+1$ equations:
$$\sum p_i = (d+1)vol(CH(S)), \qquad \sum x_ip_i = (d+1)x_c vol(CH(S)),$$
where $x_c$ is center (center of mass) of $CH(S).$ Since the right hand sides
of these equations not depend of triangulation that in fact
$p(t) \in {\bf R}^{n-d-1}.$  Let

$$Q = \mbox{convex hull}\{p(t): t \mbox{\quad is a triangulation of}
\quad S\},$$
i.e. $Q$ is a convex hull of the set of images in ${\bf R}^n$ all
triangulations of $S.$ The convex polytope $Q$ called {\em secondary
polytope} and denoted by $\sum(S).$ The dimension of this polytope is
$n-d-1.$

Even for a simple configuration of $S$ the secondary polytope is not simple.
When $S$ in the plane consist of $n$ vertices of a convex $n$--gon,
then secondary polytope called {\em associahedron} [12]. Associahedron is
a simple (n-3)-dimensional polytope with $n(n-3)/2$ facets. For example,
associahedron of a pentagon is a pentagon, and
associahedron of a hexagon is a simple 3-polytope with 14 vertices, 21 edges,
and 9 facets. An interesting application of the associahedron to the
theoretical Computer Science has been given by Sleator, Tarjan, and
Thurston [20]. These authors derive a tight upper bound for the rotation
distance between binary trees with $n$ nodes by proving that the diameter of
the associahedron equals $2n-10,$ for large $n.$

A triangulation $t$ of $S$ is said to be {\em regular} if there exists
a function on $P$ that is piecewise linear and strictly convex with respect
to $t$. (A convex piecewise linear function over a triangulation $t$ is
said to be {\em strictly convex} if it is given by a different linear
function on each maximal cell of $t$).

Regular triangulations are really just the duals of {\em power diagrams},
and Edelsbrunner call them {\em weighted} Delaunay triangulations [8].

There are several equivalent ways to define the notion of a regular
triangulation of $S$. For example, Gel'fand, Kapranov, and Zelevinski [5]
call regular a following triangulation $t$:

Choose numbers $y_1, \ldots, y_n$ and let $W = \{(x_1,y_1), \ldots,
(x_n, y_n)\} \subset {\bf R}^{d+1}$. If $\bar \Delta=\{(x_{i_1}, y_{i_1}),
\ldots, (x_{i_{d+1}}, y_{i_{d+1}})\}$ is a facet of $W$ in the
{\em lower hull} of $W$ (i.e., the last component of the outward normal of
its supporting hyperplane is negative) then $\Delta=\{x_{i_1},
\ldots, x_{i_{d+1}}\}$ is a $d$--face of the triangulation $t$.

Let $y_i = ||x_i||^2.$ Then we get the Delaunay triangulation of $S$.
Therefore, from this definition the Delaunay triangulation is regular.

Consider a set $S$ of $d+2$ points in ${\bf R}^d$.
From Gale diagram [7] follows (see also Schlegel [18], Lawson [11])
that there are exactly two ways to triangulate $S$.
Indeed, the two ways correspond to the two sides (lower and upper) of the
$(d+1)$--simplex that is convex hull of corresponding lifted points in
${\bf R}^{d+1}.$ A {\em flip} is the operation that substitutes one
triangulation of $S$ for the other [8].

There are three types of flips in two dimensions, and we denote a flip
by the number of triangles before and after the flip. So the flips in two
dimension are of type '1 to 3', '2 to 2', and '3 to 1'. The first type
introduces a new point, and the last type removes a point. The flips in
three dimensions can be classified as '4 to 1', '3 to 2', '2 to 3', and
'1 to 4'.

Now we can give the main results of this section. First theorem belongs to
Gel'fand, Kapranov, and Zelevinsky. They used another terminology.

\begin{theorem} The vertices of the secondary polytope $Q = \sum(S)$ are in
one-to-one correspondence with the regular triangulations of $S$, and
the edges of $Q$ are corresponding to flips.
\end{theorem}

    From this theorem and Section 1 follow:

\begin{theorem} There is a sequence of flips that connect any regular
triangulation of $S$ with the Delaunay triangulation of $S$, and after each
flip decreases parabolic functional.
\end{theorem}

\begin{proof} It is easy to see that for any triangulation $t$ a parabolic
functional $Vr(t)=\sum ||x_i||^2 p_i$, where $p_i = \sum vol (\Delta_{ij})$
as above. $p_i$ is $i$--th coordinate of ${\bf R}^n$ and therefore
$Vr$--functional is a linear function on ${\bf R}^n$. Secondary polytope
$Q = \sum(S)$ is a convex polytope in ${\bf R}^n$, and Delaunay triangulation
of $S$ ($DT(S)$) is a vertex of $Q$ that gives minimum for function $Vr$.
It is clear how to find sequence of neighboring vertices of $Q$ decreasing
parabolic functional and connected any regular triangulation (vertex of $Q$)
with vertex corresponding to the Delaunay triangulation.
\end{proof}

\section {Incremental construction of the Voronoi diagram.}

A natural idea is to construct the Voronoi diagram
by {\em incremental insertion}, i.e. to obtain $V(S)$ from
$V(S)\backslash\{x\}$
by inserting the point $x$. The insertion process is, maybe,
better described, and implemented in the dual environment, for the Delaunay
triangulation: construct $DT_i = DT(\{x_1, x_2,..., x_{i-1}, x_i\})$ by
inserting the point $x_i$ into $DT_{i-1}$. The advantage over a direct
construction of $V(S)$ is that Voronoi vertices that appear in intermediate
diagrams but not in the final one need not be constructed and stored.

    Several algorithms proposed for Delaunay triangulation are based on the
notion of a local transformation henceforth referred to as a flip.
Historically the first such algorithm is due to Lawson [10]. Given a finite
point set in the plane, the algorithm first construct an arbitrary
triangulation of the set. This triangulation is then gradually altered
through a sequence of edge-flips until the Delaunay triangulation is
obtained. The generalization of this method to  ${\bf R}^d$, $d>2$ has
difficulties, and it is incorrect if the flips are applied to an arbitrary
initial triangulation [8]. Joe [9] shows that if a single point, $x_i$, is
added to the Delaunay triangulation $DT_{i-1}$ in ${\bf R}^3$ then many
different sequences of flips will succeed in constructing the Delaunay
triangulation $DT_i$. This can be used as the basis of an incremental
algorithm. Rajan [16] considers Delaunay triangulation in arbitrary
dimensions, ${\bf R}^d$, and argues that a single point can always be
added by a sequence of flips. However, he needs a priority queue to find
the appropriate sequence, which takes logarithmic time per flip.
Edelsbrunner and Shah [8] using "weighted points" method unifies and
extends the algorithmic results of Joe [9] and Rajan [16]. In particular,
they show that many different sequences of flips can be used to add a single
point to a regular triangulation in ${\bf R}^d$. This eliminates the need
for a priority queue that sorts the flips. This section we show how to
applied parabolic functional and secondary polytope for construction of
the Voronoi diagram incrementally.

The algorithm constructs the Delaunay triangulation of a give set
$S = \{x_1, x_2, ..., x_n\}$ in ${\bf R}^d$ incrementally. It is convenient
to first construct an artificial $d$--simplex
$S_0 = \{x_d, x_{-d+1},...,x_0\}$, so that $S$ is contained in it. The
$d+1$ artificial points can be conveniently chosen at infinity, so that
choice of points guarantees that $DT(S)$ is a subcomplex of
$DT(S \bigcup S_0)$ [8]. In fact, $DT(S)$ consist of all simplices of
$DT(S \bigcup S_0)$ that are not incident to any point of $S_0$.

Let $t$ be a triangulation of $S$. Call $t$ {\em locally non--optimal}
triangulation if there is a convex subcomplex $\sigma \in t$ that after
the flip $\sigma$ a new triangulation is {\bf regular} and {\bf decreases}
parabolic functional (see section 1.) i.e. $Vr(t) > Vr(t')$, where $t'$
is the new triangulation. In accordance with this definition we call
$t$ L.O.T. (locally optimal triangulation) if $t$ is not locally non--optimal.
We denote a locally optimal triangulation of $S$ as $LOT(S)$. We will show
later that $LOT(S_i) = DT_i$ for all $i=0,1,...,n$.

\medskip

{\em THE INCREMENTAL ALGORITHM FOR CONSTRUCTION OF} $LOT(S)$

1 Construct $LOT(S_0)$;

2 {\bf for} $i := 1$ {\bf to} $n$ {\bf do}

3  locate the d-simplex $s$ in $LOT(S_{i-1})$ that contains $x_i$

4 {\bf if} $LOT(S_{i-1} \bigcup \{x_i\}$ is not L.O.T. {\bf then}

5 flip $LOT(S_{i-1} \bigcup \{x_i\}$;

6 {\bf while} there exist non--optimal $d+2$--subcomlex {\bf do}

7   find a non--optimal $d+2$--subcomplex $\sigma$;

8   flip $\sigma$

9    {\bf endwhile}

10    {\bf endif}

11 {\bf endfor}

\medskip

This algorithm could fail for two reason. First, it could be that the
{\bf while} loop does not terminate, because it cycles in an infinite loop
of flips. Second, if the algorithm would stop before reaching the $DT(S)$.
We show this cannot happen.

The proof that sequence of flips in the algorithm does not cycle easily
follows from the consideration of the parabolic functional $Vr(t)$. Indeed,
this functional decreases after each flip in the algorithm.

For each $i$ in algorithm $LOT(S_i)$ is locally optimal triangulation of
$S_i.$ It means, that $LOT(S_i)$ is regular triangulation, and there is
no a locally non-optimal subcomplex $\sigma.$ Then $LOT(S_i)$ is a vertex
of a secondary polytope $Q_i = \sum(S_i).$ $DT_i$ is vertex of $Q_i$ also,
and therefore from Theorem 4 (see section 2.) follows that there is
sequence of flips that connect $LOT(S_i)$ and $DT_i.$ From theorem 1 follows
that $DT_i$ gives minimum for parabolic functional. From other side parabolic
functional cannot be decrease for $LOT(S_i)$. Consequently $LOT(S_i)$ is
$DT_i$, and therefore

\begin{theorem} Constructed in the algorithm a locally optimal triangulation
of $S$ ($LOT(S)$) is the Delaunay triangulation of $S$.
\end{theorem}

\section {The Dirichlet functional on triangulations}

Let $S = \{x_i\}$ be a set of $n$ points in ${\bf R}^d$, each
associated with a real number $y_i$.
Denote by $Y$ the set of these numbers, i.e. $Y = (y_1,..., y_n)$.
There are a lot of different problems in Geography,
Geology, Topography, CAD/CAM etc., where we need to construct a surface in
${\bf R}^{d+1}$ corresponding to this dataset. The main problem is the
following: to find a function $y = f(x)$, such that $f(x_i) = y_i$. One of
the oldest and the most famous methods is modeling by triangulation. If we
have some triangulation of $S$ then for a set of data $Y$ there is only one
method to construct a piecewise linear function (polyhedral surface)
on this triangulation. Usually Delaunay triangulation is used for
this purpose.

One of the minimum criterion is a discrete analogue of the Dirichlet
functional:
$\int ||grad\, f(x)||^2 dx$. For interval ($d=1$), a spline of $\deg=2k-1$ is
a function $y=f(x)$ such that $f(x_i)=y_i$ and
$$\int\limits_a^b [f^{(k)}(x)]^2dx = min.$$
For $k=1, d>1$ and piecewise linear function $f$ we get
$$DF(t,Y) = \int\limits_{CH(S)} ||grad\,f(x)||^2dx 
=\sum\limits_i \frac{(vol(\Delta_i(Y)))^2}{vol(\Delta_i)} - vol(CH(S)).$$
The triangulation of $S$ that minimizing the functional $DF$ can be  called the
{\ em discrete spline triangulation} (DST).
For the plane DST does not depend upon $Y$ and it is DT. Rippa [11] for
$d=2$ proved that $DF(t,Y)$ achieves its minimum iff $t$ is DT.

For $d>2$ the Delaunay triangulation of $S$ could be not optimal for this
functional. Let $$x_0=(0,0,\cdots,0);\quad
x_i=(\underbrace{0,\cdots,0}_{i-1}, 1, \underbrace{0,\cdots,0}_{n-i}),
i=1,\cdots,d; \quad x_{d+1}=(a,\cdots,a).$$
There are exactly two distinct triangulations of $S=(x_0,x_1,\cdots,x_{d+1}).$
Denote by $t_1$ a triangulation consist of two simplices:
$\Delta_1=(x_0, x_1,\cdots,x_d); \Delta_2=(x_{d+1}, x_1,\cdots,x_d),$ and
by $t_2$ another one. Then $t_1$ is DT iff $a>1$. It is easy to show by
direct calculation that $t_1$ is DST iff $a>\frac{1}{d-1}.$ Therefore,
for $\frac{1}{d-1}<a<1\; t_1$ is DST, but is not DT.

The proof of  Rippa's theorem directly follows from the fact that
the $DF$--functional for triangulation $t$ of a quadrilateral is minimum
if $t$ is DT. The proof also follows from some general result that is
given below.

Let $S$ be a set of $d+2$ points $x_1,..., x_{d+2}$ in ${\bf R}^d$.
Suppose $S$ admits two triangulations $t_1$ and $t_2$,
and $Y = (y_1,..., y_{d+2})$ is a set numbers corresponding to
$x_1,..., x_{d+2}$ as above. Let $B(Y,S) = DF(t_1,Y)-DF(t_2,Y)$. Note
$B(Y,S)$ is a quadratic form depending on $Y$.

\begin{theorem} The optimal (DST) triangulation of $S$ for $n=d+2$ does not
depend on $Y$.
\end{theorem}

\begin{proof}
It is easy to see that for arbitrary set of real numbers
$(a_0,a_1,\cdots,a_d):$

$$B(\hat Y,S)=B(Y,S),\quad\mbox{where}\quad \hat Y=(\hat y_1,\cdots,\hat y_n),
\quad\mbox{and}$$
$$\hat y_i=y_i+a_0+\sum\limits_{j=1}^d a_jx_{ij}; \quad
x_i=(x_{i1},\cdots,x_{id}) \in S.$$

Then for $d+1$--dimensional subspace  $\Re \subset {\bf R}^n$ that is the 
linear hull of the set of vectors:
$$(1,\cdots,1), \quad (x_{11},\cdots,x_{n1}),\cdots, (x_{1d},\cdots,x_{nd})$$
the quadratic form $B(Y,S)$ is vanished, and $B(Y+Z,S)=B(Y,S)$ if
$Z\in\Re.$ Note that $n=d+2$ and $\dim{\Re}=d+1$ therefore $B(Y,S)$ could be not vanish
only on 1--dimensional subspace of ${\bf R}^n$ that is orthogonal to $\Re.$
Thus
$B(Y,S) = const(S)L^2(Y),$ where $L(Y)$ is some linear form on
$Y\in{\bf R}^n$. Therefore sign of the $B(Y,S)=DF(t_1,Y)-DF(t_2,Y)$ does
not depend on $Y$, and if $t_1$ is optimal for some $Y$ (i.e.
$DF(t_1,Y)<DF(t_2,Y)$) that it is optimal for any $Y.$
\end{proof}

We state an open problem concerning DST:

{\em -- Does DST depend on $Y,$ when $n>d+2$ ?}

\section {Concluding Remarks.}

    Voronoi diagram and Delaunay triangulation have a fair number of
applications, including the generation of grids for point configuration and
for surface interpolation. No doubt that in two-dimension these tesselations
are the best for these purposes. Indeed, the main motivation for studying
the problems solved in this paper is our intention to implement Voronoi and
Dirichlet functionals in dimensions beyond ${\bf R}^3.$ The optimal
triangulations for these functionals in ${\bf R}^d, d>2$ could be not the
equal.  We do not know is optimal triangulation for Dirichlet functional is
regular?
If it is regular then for it construction the algorithm in section 3 is
suitable. Instead of $Vr$--functional there have to be used $DF$--functional.
It would be interesting to study dual tessellation for DST.
In other word, what is analogue of Voronoi diagram for DST?

\medskip

\medskip

{\bf\large Acknowledgement.} The author wish to thank Herbert Edelsbrunner for
helpful discussions concerning this paper. Especially, I am grateful to him
for explaining his and Shah work [8] to me. It is easy to find results of
these explanations and influence of the paper [8] in section 3.

\medskip

\medskip

{\Large\bf References}
\medskip

[1] Aurenhammer F., Klein R. Voronoi Diagrams // Optimierung und Kontrolle,
Bericht Nr. 92, Karl-Franzens-Univ. Graz and Tech. Univ. Gratz, 1996.

[2] D'Azevedo, E.F. and Simpson, R.B.
On optimal interpolation triangle incidences, SIAM J. Sci. Statist. Comput.,
vol. 10, No. 6 , pp. 1063--1075, 1989.

[3] D'Azevedo, E.F. Optimal triangular mesh generation by coordinate
transformation, SIAM J. Sci. Statist. Comput., vol. 12 , No. 4,
pp. 755--786, 1991.

[4] Billera L.J., Filliman P., Sturmfels B. Construction and complexity
of secondary polytopes. Advances in Math., vol. 83, No. 2, pp. 155-179, 1990.

[5] Gel'fand I.M., Kapranov M.M., Zelevinsky A.V. Newton polyhedra of
principal {\em A}--determinants, Soviet Math. Dokl. {\bf 308}, pp. 20-23,
1989.

[6] Delaunay, B.N. Sur la sph{\`e}re vide. A la memoire de Georges Voronoi.
Izv. Akad. Nauk SSSR, Otd. Mat. i Estestv. nauk, No 7, pp. 793-800, 1934.

[7] Emelichev, V.A., Kovalev, M.M., Kravtsov, M.K. Polytopes, Graphs,
Optimization. Moscow: "Nauka", 344 p., 1981 (in Russian).

[8] Edelsbrunner, H., Shah, N.R. Incremental topological flipping works
for regular triangulation. Algorithmica, 15, pp. 223-241, 1996.

[9] Joe, B. Construction of three -- dimensional Delaunay triangulation using
local transformation, Computer Aided Geometric Design, No. 8, pp. 123-142,
1992.

[10] Lawson, C.L. Software for $C^1$ surface interpolation, in: J.R. Rice, ed.,
Mathematical Software III, Academic Press, New York, pp. 161-194, 1977.

[11] Lawson, C.L. Properties of n-dimensional triangulations,
Computer Aided Geometric Design, No. 3, pp. 231-246, 1986.

[12] Lee, C. The associahedron and triangulations of the $n$--gon.
European J. Combin., vol. 10, pp. 551-560, 1989.

[13] Musin, O.R. Delaunay triangulation and optimality, ARO Workshop
Comp. Geom., Raleigh (North Carolina), pp. 37-38, 1993.

[14] Musin, O.R. Index of harmony and Delaunay triangulation,
Symmetry: Culture and Science, Vol. 6, No. 3, pp. 389-392, 1995.

[15] Musin, O.R. Properties of the Delaunay triangulation,
Proc. 13th Annu. ACM Sympos. Comput. Geom., pp. 424-426, 1997.

[16] Rajan V.T. Optimality of the Delaunay triangulation in {$R^{d}$},
Proc. 7th Annu. ACM Sympos. Comput. Geom., pp. 357--363, 1991.

[17] Rippa, S. Minimal roughness property of the Delaunay triangulation,
Computer Aided Geometric Design, No. 7, pp. 489-497, 1990

[18] Schlegel, V. Ueber die verschiedenen Formen von Gruppen, welcher
beliebige Punkte im n--dimensionalen Raum bilden konnen. -
Arch. Math. Phys., 10, 1891.

[19] Sibson, R. Locally eguiangular triangulations, Comput. J.,
Vol.21, No. 3, pp. 243-245, 1978.

[20] Sleator, D.D., Tarjan, R.E., Thurston, W.P. Rotation distance,
triangulations, and hyperbolic geometry, J. Amer. Math. Soc., vol. 1,
pp. 647-681, 1988.

[21] Voronoi, G.F. Nouvelles applications des param{\`e}tres continus
{\`a} la th{\'e}orie des formes quadratiques, J. Reine u. Angew. Math.,
34 , 198-287, 1908.
\end{document}